\input ./style/arxiv-vmsta.cfg
\documentclass[numbers,compress,v1.0.1]{vmsta}

\volume{4}% Updated by VTEXPTS2LaTeX.exe, 06.04.2017 12:49
\issue{1}% Updated by VTEXPTS2LaTeX.exe, 06.04.2017 12:49
\pubyear{2017}
\firstpage{79}% Updated by VTEXPTS2LaTeX.exe, 06.04.2017 12:49
\lastpage{89}% Updated by VTEXPTS2LaTeX.exe, 06.04.2017 12:49
\doi{10.15559/17-VMSTA75}% Updated by VTEXPTS2LaTeX.exe, 04.04.2017
\newtheorem{Th}{Theorem}
\newtheorem{Th_2}[Th]{Theorem}

\theoremstyle{definition}
\newtheorem{rem}{Remark}

\startlocaldefs
\allowdisplaybreaks

\def\rd{\rm d}
\def\bZ{\mathbf Z}

\def\bR{\bf R}

\def\bN{\bf N}
\def\hX{\hat X}

\def\hY{\hat Y}
\def\be{\bar e}
\def\tr{\tilde r}

\endlocaldefs

\begin{document}
\begin{frontmatter}

\title{Asymptotic behaviour of non-isotropic random walks with heavy tails}

\author[a]{\inits{M.}\fnm{Mark}\snm{Kelbert}\corref{cor1}}\email
{mkelbert@hse.ru}
\cortext[cor1]{Corresponding author.}
\address[a]{National Research University Higher School of Economics,
Moscow RF}

\author[b]{\inits{E.}\fnm{Enzo}\snm{Orsingher}}\email
{enzo.orsingher@uniroma1.it}
\address[b]{Sapienza University of Rome}

\markboth{M. Kelbert, E. Orsingher}{Asymptotic behaviour of
non-isotropic random walks with heavy tails}

\begin{abstract}
A random flight on a plane with non-isotropic displacements at the
moments of direction changes is considered. In the case of
exponentially distributed flight lengths a Gaussian limit theorem is
proved for the position of a particle in the scheme of series when jump
lengths and non-isotropic displacements tend to zero. If the flight
lengths have a folded Cauchy distribution the limiting distribution of
the particle position is a convolution of the circular bivariate Cauchy
distribution with a Gaussian law.
\end{abstract}

%
%\begin{keywords}[2010]% [PACS], [JEL]
%\kwd{}
%\kwd{}
%\kwd{}
%\kwd{}
%\end{keywords}

%
\received{24 November 2016}% Updated by VTEXPTS2LaTeX.exe, 04.04.2017
%09:30
%
\revised{9 March 2017}% Updated by VTEXPTS2LaTeX.exe, 04.04.2017 09:30
\accepted{14 March 2017}% Updated by VTEXPTS2LaTeX.exe, 04.04.2017
%09:30
\publishedonline{6 April 2017}
\end{frontmatter}

\section{Introduction}\label{Intro}

We consider the problem of random flights in Euclidean spaces defined
by a series of displacements, $\bar r_j$, the magnitude and direction of
each one being independent of all the previous ones. This model was
introduced by Karl Pearson in 1905 and has a long and interesting
history, both as a purely mathematical problem in probability
theory and as a model for various physical and chemical processes \cite
{C}. The majority of papers %tackle
investigate the problem of random flights with %distribution %of
the orientation of movements uniformly
distributed over a sphere, and deviations separated by exponentially
distributed or Dirichlet distributed time lapses (cf. discussions in
\cite{DGO,LC,ODG}).
For the most recent developments in the studying of the random walks in
a random environment we refer to \cite{DK} and the papers cited therein.
In this short note
we introduce a novel feature in the form of non-isotropic displacements
at the moments of the direction changes. As a model of this
non-isotropic perturbation we consider Hadamard's (or componentwise) product
of a fixed deterministic vector $(\varDelta_1,\varDelta_2)$ with the unit
vector $\be=(\cos\theta_j,\sin\theta_j)$ in the direction of the
previous movement. In the following analysis
these perturbations are assumed to be small and the direction changes
are frequent enough. In this note we are mainly interested in the case
where the distribution of the i.i.d. flight lengths has a heavy tail,
say it follows a folded Cauchy distribution.

More precisely, we consider a planar, non-isotropic random walk
performed by a particle taking steps $(X_j, Y_j), j\in\bN$. We assume that
\begin{equation}
X_j=(R_j+\varDelta_1)\cos\theta_j, \qquad Y_j=(R_j+\varDelta_2)\sin\theta_j, \label{1}
\end{equation}
where $\theta_j$ and $R_j=|\bar r_j|$ are independent positive random
variables (hereafter r.v.'s), and $\varDelta_1\neq\varDelta_2$ are
deterministic positive real numbers. We assume that
$\theta_j$ are uniformly distributed in $[0,2\pi) $ and positive
r.v.'s $R_j$ are identically distributed with density $f(r), r>0$.
Clearly, after $n$ steps the position reached by the moving particle is
given by
\begin{equation}
\hX_n=\sum\limits_{j=1}^n X_j, \qquad \hY_n=\sum\limits_{j=1}^n Y_j. \label{2}
\end{equation}
A possible sample path of the random walk (\ref{2}) is depicted in
\begin{figure}[h!]
\includegraphics{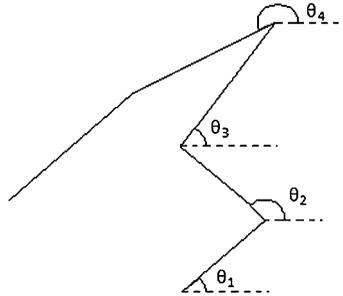}
\caption{A sample path of the random walk}\label{Figure 1}
\end{figure}
Figure 1 and can be interpreted as the position of a particle taking
jumps at integer-valued times, with arbitrary orientation.
Since $\varDelta_1\neq\varDelta_2$,
the distribution of the random walk $(\hX_n,\hY_n), n\geq1$,
as well as that of its asymptotic limiting process,
is not rotation invariant.
If the angles $\theta_j$ are non-uniformly distributed on $[0,2\pi)$
the resulting random motion is anisotropic as well, this case will be
studied elsewhere.
Here two qualitatively different examples are considered:
1) the exponential distribution of the i.i.d. flight lengths, and 2)
the folded Cauchy distribution when all the moments $m_r, r\geq1$ are infinite.
Naturally, in the former case under a suitable
scaling one obtains the Gaussian limit with independent components
having different variances. In the latter case the limiting law is a
convolution of
a circular bivariate Cauchy distribution with a Gaussian law.

\section{Main results}\label{Results}

In case 1)
we are working under the following assumptions
\begin{itemize}
\item(i) The jump lengths $R_j=|\bar r_j| , j=1,\ldots, n$ are
exponentially distributed with parameter
\begin{equation}
\mu^{(n)}=\frac{\mu}{t} n^{-1/2}. \label{3}
\end{equation}
\noindent
\item(ii) The asymmetry conditions: for some $C_1,C_2>0$ (we are
interested in the case $C_1\neq C_2$)
\begin{equation}
\varDelta_i^{(n)}=C_i n^{-1/2},\quad i=1,2, \label{4}
\end{equation}
\end{itemize}
i.e. the displacement vectors decrease with $n$ and are the same for
all $1\leq j\leq n$.
Condition (i) means that for fixed values of $n$ the step lengths $R_j$
are i.i.d. with exponential distribution whose parameter $\mu^{(n)}$
is adjusted continuously.
One can easily see that ${\bf E} X_j={\bf E} Y_j=0$. Next,
\begin{equation}
{\bf E}X_1^2= \bigl({\bf E}R_1^2+2
\varDelta_1 {\bf E}R_1+\varDelta_1^2
\bigr){\bf E}\bigl[(\cos\theta)^2\bigr]\approx n^{-1}
\biggl(\frac{t^2}{\mu^2}+\frac
{C_1t}{\mu}+\frac{C_1^2}{2} \biggr), \label{5a}
\end{equation}
in view of the equality ${\bf E}[(\cos\theta)^2]=\frac{1}{2}$. In
fact, (\ref{5a}) is an identity.
%below
Here and in what follows the symbol $\approx$ is used
to indicate that LHS is an expansion of RHS up to $O(n^{-k})$ for some
$k\in\bZ$ but the difference RHS--LHS is $o(n^{-k})$. Similarly,
\[
{\bf E}Y_1^2\approx n^{-1} \biggl(
\frac{t^2}{\mu^2}+\frac{C_2t}{\mu
}+\frac{C_2^2}{2} \biggr).
\]
Moreover, $X_j$ and $Y_j$ are dependent but not correlated r.v.'s.
These facts suggest that a joint limiting distribution is Gaussian and is
represented in the form of two independent diffusions. The proof may be
provided by the standard methods via the checking
%of
Lindeberg's conditions.
However, we prefer to use a direct computation to pave the way for
further results for r.v.'s with the heavy tails.
\begin{Th}
Under the assumptions (\ref{3}) and (\ref{4}) the sequence $(\hX
_n,\hY_n)$ defined in (\ref{2}) weakly converges to the zero-mean
Gaussian vector
$(X,Y)$ where $X, Y$ are independent and possess the variances
\begin{equation}
{\rm Var}(X)=\frac{t^2}{\mu^2}+\frac{C_1t}{\mu}+\frac{C_1^2}{2}, \qquad {\rm
Var}(Y)=\frac{t^2}{\mu^2}+\frac{C_2t}{\mu}+\frac{C_2^2}{2}. \label{5}
\end{equation}
\end{Th}
Our main result in Theorem 2 below establishes the limiting law for the
(folded) Cauchy flights:
\begin{equation}
f(r)=\frac{2}{\pi}\frac{a}{r^2+a^2},\quad  a>0, r>0. \label{6a}
\end{equation}
Let us remind that the standard circular bivariate Cauchy distribution
has the joint PDF (see \cite{F}, Ch.II, formula (1.19) or \cite{Fe})
\begin{equation}
f(x_1,x_2)=\frac{1}{2\pi}\frac{1}{(1+x_1^2+x_2^2)^{3/2}}.
\label{4a}
\end{equation}
\begingroup
\abovedisplayskip=7.5pt
\belowdisplayskip=7.5pt
\begin{Th_2}
Assume the condition (\ref{4}) and fix $b>0$. Let the parameter of the
folded Cauchy distribution (\ref{6a}) for the jump lengths be scaled
as $a_n=\frac{\pi b}{2n}$. Then
the distribution of random vector $(\hX_n,\hY_n)$ weakly converges as
$n\to\infty$ to the convolution of the cumulative distribution functions
$F_{X,Y}\circ F_{V,W}$ where $(X,Y)$ is a zero mean Gaussian vector
with independent components,
\begin{equation}
{\rm Var}(X)=\frac{C_1^2}{2},\qquad  {\rm Var}(Y)=\frac{C_2^2}{2}, \label{6b}
\end{equation}
and the vector $(V,W)$ has a circular bivariate Cauchy distribution
with the shape parameter $b$, i.e.
\begin{equation}
f_{V,W}(x_1,x_2)=\frac{1}{2\pi}
\frac{b}{(b^2+x_1^2+x_2^2)^{3/2}}. \label{4b}
\end{equation}
\end{Th_2}
\begin{rem}
Let us represent the position of the moving particle as
a result of $n$ random flights and $n$ non-isotropic displacements
$\hX_n=U_n+T_n, \hY_n=V_n+S_n$. Here $U_n=\sum_{j=1}^nR_j\cos\theta_j, T_n=\varDelta_1\sum_{j=1}^n\cos
\theta_j,
V_n=\sum_{j=1}^nR_j\sin\theta_j, S_n=\varDelta_2\sum_{j=1}^n\sin\theta_j$. Then
\begin{equation}
{\bf E} e^{i\alpha\hX_n+i\beta\hY_n}={\bf E} \bigl[e^{i(\alpha
U_n+\beta V_n)}e^{i(\alpha T_n+\beta S_n)} \bigr].
\label{4c}
\end{equation}
The expectations in the RHS of (\ref{4c}) may be split as $n\to\infty
$, cf. (\ref{21}) below. Moreover $T_n$ and $S_n$ are asymptotically
independent.
However, the pair $(U_n,V_n)$ asymptotically follows a circular
bivariate Cauchy law.
\end{rem}

%% Theorem %%
%%%%%%%%%%%%%%%%%%%%%%
%%\begin{thm}[]\label{} Theorem
%%\end{thm}

\section{Proofs}\label{Proofs}
The initial steps are the same for both Theorems 1 and 2 and valid for
any PDF $f(r)$ of %the
flight lengths.
They are based on the properties of Bessel functions $J_{\nu}(x)$
which are solutions of ODEs
\begin{equation}
Ly=x^2\frac{{\rd}^2y}{{\rd} x^2}+x\frac{{\rd} y}{{\rd}
x}+\bigl(x^2-
\nu^2\bigr)y=0,
\end{equation}
and admit the expansion \cite{AS,W}
\begin{equation}
J_{\nu}(x)=\sum\limits
_{m=0}^{\infty}
\frac{(-1)^m}{m!\varGamma
(m+k+1)} \biggl(\frac{x}{2} \biggr)^{2m+\nu}.
\label{12a}
\end{equation}
Let us fix a small open neighbourhood $U$ of $(0,0)$. For $(\alpha
,\beta)\in U$
the characteristic function of the steps $(X_j,Y_j)$ reads
\begin{align}
\varphi(\alpha,\beta)&={\bf E} \bigl[\exp{ \bigl(i\alpha[R+\varDelta _1]\cos\theta+i\beta[R+\varDelta_2]\sin\theta \bigr)} \bigr]\nonumber\\
&=\frac{1}{2\pi}\int_0^{2\pi}{\rm d}\theta\int_0^{\infty}e^{i\alpha(r+\varDelta_1)\cos\theta+i\beta(r+\varDelta_2)\sin\theta}f(r)\,{\rm d}r\nonumber\\
&=\int_0^{\infty}J_0 \Bigl(\sqrt{\bigl(\alpha^2+\beta^2\bigr)r^2+2r\bigl(\alpha ^2\varDelta_1+\beta^2\varDelta_2\bigr)+\bigl(\alpha^2\varDelta_1^2+\beta^2\varDelta _2^2\bigr)} \Bigr)f(r)\,{\rm d}r.
\label{6}
\end{align}
\endgroup

Due to the addition formula of Bessel functions (\cite{GR}, formula
8.531, page 979)
\begin{equation}
J_0\bigl(\sqrt{\tr^2+\rho^2-2\tr\rho\cos
\phi}\bigr)=J_0(\tr)J_0(\rho )+2\sum
\limits_{k=1}^{\infty}J_k(
\tr)J_k(\rho)\cos(k\phi), \label{7}
\end{equation}
where, in our case,
\begin{eqnarray}
&\displaystyle\tr^2=\bigl(\alpha^2+\beta^2\bigr)r^2, \qquad \rho^2=\bigl(\alpha^2\varDelta_1^2+\beta ^2\varDelta_2^2\bigr), &\nonumber\\
&\displaystyle\cos\phi=-\frac{\alpha^2\varDelta_1+\beta^2\varDelta_2}{\sqrt{\alpha^2+\beta^2}\sqrt{\alpha^2\varDelta_1^2+\beta^2\varDelta_2^2}}. \label{8}&
\end{eqnarray}

\section{Proof of Theorem 1}\label{Theorem 1}
For the sake of brevity we omit the upper index $\mu^{(n)}$, low index
$\varphi_n(\alpha,\beta)$, etc., whenever it is possible.
We can %work out
calculate explicitly the characteristic function $\varphi(\alpha
,\beta)$ by means of
%an
integration term by term. Further, we must keep into account the
additional result
\begin{equation}
\int_0^{\infty}e^{-\alpha x}J_{\nu}(
\beta x)\,{\rm d}x=\frac{
[\sqrt{\alpha^2+\beta^2}-\alpha ]^{\nu}}{\beta^{\nu}\sqrt
{\alpha^2+\beta^2}}, \quad \nu>-1, \alpha>0 \label{9}
\end{equation}
(\cite{GR}, formula 6.611, page 707).

In view of all these formulas we have that for $f(r)=\mu e^{-\mu r}$
\begin{align}
&\varphi(\alpha,\beta)\nonumber\\
&\quad =J_0 \Bigl(\sqrt{\alpha^2\varDelta_1^2+\beta^2\varDelta_2^2} \Bigr)\int_0^{\infty}\mu e^{-\mu r}J_0\bigl(r\sqrt{\alpha ^2+\beta^2}\bigr)\,{\rm d}r\nonumber\\
&\qquad +2\sum\limits_{k=1}^{\infty}J_k\Bigl(\sqrt{\alpha^2\varDelta_1^2+\beta^2\varDelta_2^2}\Bigr)\cos \biggl(k\; {\rm arcos} \biggl[-\frac{\alpha^2 \varDelta_1+\beta^2\varDelta_2}{\sqrt{\alpha^2+\beta^2}\sqrt{\alpha^2\varDelta_1^2+\beta^2\varDelta_2^2}} \biggr] \biggr)\nonumber\\
&\qquad \times\int_0^{\infty}J_k\bigl(r\sqrt{\alpha^2+\beta^2}\bigr)\mu e^{-\mu r}\,{\rm d}r\nonumber\\
&\quad =J_0 \Bigl(\sqrt{\alpha^2\varDelta_1^2+\beta^2\varDelta_2^2} \Bigr)\frac{\mu}{\sqrt{\mu^2+\alpha^2+\beta^2}}\nonumber\\
&\qquad +2\mu\sum\limits_{k=1}^{\infty} J_k \Bigl(\sqrt{\alpha^2\varDelta _1^2+\beta^2\varDelta_2^2} \Bigr)\cos \biggl(k\; {\rm arcos} \biggl[-\frac{\varDelta_1\alpha^2+\varDelta_2\beta^2}{\sqrt{\alpha^2+\beta^2}\sqrt{\alpha^2\varDelta_1^2+\beta^2\varDelta_2^2}} \biggr] \biggr)\nonumber\\
&\qquad \times\frac{(\alpha^2+\beta^2)^{-k/2}}{\sqrt{\alpha^2+\beta^2+\mu^2}} \bigl(\sqrt{\mu^2+\alpha^2+\beta^2}-\mu \bigr)^k.\label{10}
\end{align}

A crucial point is now to preserve only the relevant terms of the
expansion of $\varphi_n(\alpha,\beta)$ in view of the evaluation of
the limit for the
characteristic function\break
$\lim_{n\to\infty}[\varphi_n(\alpha,\beta)]^n$, taking into
account that $\mu=\mu^{(n)}, \varDelta_i=\varDelta_i^{(n)}, i=1,2$. Since
\begin{equation}
\frac{ (\sqrt{\mu^2+\alpha^2+\beta^2}-\mu )^k}{\sqrt{\mu
^2+\alpha^2+\beta^2}}= \bigl(\sqrt{\mu^2+\alpha^2+
\beta^2} \bigr)^{k-1} \biggl(1-\frac{\mu}{\sqrt{\mu^2+\alpha^2+\beta^2}}
\biggr)^k,
\end{equation}
we can cut the terms of the expansion for $k\geq2$. To justify this
fact let us use the expansion (\ref{12a}).
All the terms with $k\geq2$ in (\ref{10}) contain the factors
$(\alpha^2\varDelta_1^2+\beta^2\varDelta_2^2)^{k/2}$ and in view of assumptions
(\ref{3}) and (\ref{4}) it is easy to check that the
sum of these terms is $o (n^{-1} )$ uniformly over $(\alpha
,\beta)\in U$. Hence,
\begin{align}
\bigl[\varphi_n(\alpha,\beta)\bigr]^n&\approx \biggl[\frac{\mu}{\sqrt{\alpha^2+\beta^2+\mu^2}}J_0 \Bigl(\sqrt{\alpha^2\varDelta_1^2+\beta^2\varDelta_2^2} \Bigr)\nonumber\\
&\quad -2J_1 \Bigl(\sqrt{\alpha^2\varDelta_1^2+\beta^2\varDelta_2^2} \Bigr)\frac{\alpha^2\varDelta_1+\beta^2\varDelta_2}{\sqrt{\alpha^2+\beta^2}\sqrt{\alpha^2\varDelta_1^2+\beta^2\varDelta_2^2}}\nonumber\\
&\quad \times\frac{\mu}{\sqrt{\alpha^2+\beta^2}}\frac{1}{\sqrt{\alpha^2+\beta^2+\mu^2}}\bigl(\sqrt{\alpha^2+\beta^2+\mu^2}-\mu\bigr) \biggr]^n\nonumber\\
&= \biggl[\frac{J_0(\sqrt{\alpha^2\varDelta_1^2+\beta^2\varDelta_2^2})}{\sqrt{1+\frac{\alpha^2+\beta^2}{\mu^2}}}\nonumber\\
&\quad -2J_1\Bigl(\sqrt {\alpha^2\varDelta_1^2+\beta^2\varDelta_2^2}\Bigr)\frac{(\alpha^2\varDelta_1+\beta^2\varDelta_2)}{(\alpha^2+\beta^2)\sqrt{\alpha^2\varDelta_1^2+\beta^2\varDelta_2^2}}\nonumber\\
&\quad \times\mu \biggl(1-\frac{\mu}{\sqrt{\alpha^2+\beta^2+\mu^2}} \biggr) \biggr]^n. \label{11}
\end{align}
For $n\to\infty$, $\varDelta_i^{(n)}\to0, i=1,2$, and since
$J_0(x)\approx1- (\frac{x}{2} )^2$ as $x\to0$ and
$J_1(x)\approx\frac{x}{2}$ for small values of $x$ we can
%extract
obtain the following relationship
\begin{align}
\bigl[\varphi_n(\alpha,\beta) \bigr]^n&\approx
\biggl[
\biggl(1-\frac{\alpha^2\varDelta_1^2+\beta^2\varDelta_2^2}{4} \biggr) \biggl(1-\frac{\alpha^2+\beta^2}{2\mu^2} \biggr)\nonumber\\
&\quad -\frac{\alpha^2\varDelta_1+\beta^2\varDelta_2}{\alpha^2+\beta^2} \mu \biggl(1-\frac{1}{\sqrt{1+\frac{\alpha^2+\beta^2}{\mu^2}}} \biggr)
\biggr]^n.
\label{12}
\end{align}
We now take the equalities
\begin{equation}
\frac{1}{\sqrt{1+x}}=\sum\limits
_{k=0}^{\infty} %
\begin{pmatrix} -1/2\cr k
\end{pmatrix} %
x^k=\sum\limits
_{k=0}^{\infty}(-1)^k
\begin{pmatrix} 2k\cr k
\end{pmatrix} %
\frac{x^k}{2^{2k}}, \label{20c}
\end{equation}
and thus for small values of $x$ we have that $(1+x)^{-1/2}\approx
1-\frac{x}{2}$.
In conclusion, by writing explicitly $\mu^{(n)}$, $\varDelta_i^{(n)},
i=1,2$, as in the assumptions (i) and (ii) we have that
\begin{align}
\bigl[\varphi_n(\alpha,\beta)\bigr]^n &\approx \biggl[ \biggl(1-\frac{\alpha^2\varDelta_1^2+\beta^2\varDelta_2^2}{4} \biggr) \biggl(1-\frac{\alpha^2+\beta^2}{2(\mu^{(n)})^2} \biggr)\nonumber\\
&\quad -\frac{\alpha^2\varDelta_1+\beta^2\varDelta_2}{\alpha^2+\beta^2}\mu^{(n)} \biggl(\frac{\alpha^2+\beta^2}{2(\mu^{(n)})^2} \biggr) \biggr]^n\nonumber\\
&= \biggl[ \biggl(1-\frac{\alpha^2\varDelta_1^2+\beta^2\varDelta_2^2}{4} \biggr) \biggl(1-\frac{\alpha^2+\beta^2}{2(\mu^{(n)})^2}\biggr)-\frac{\alpha^2\varDelta_1+\beta^2\varDelta_2}{2\mu^{(n)}} \biggr]^n\label{13}
\end{align}
(by assumptions (i) and (ii))
\begin{align}
&\approx \biggl[1-\frac{\alpha^2C_1^2+\beta^2C_2^2}{4n}-\frac{(\alpha
^2+\beta^2)t^2}{2n\mu^2}-
\biggl(\alpha^2\frac{C_1}{\sqrt{n}}+\beta ^2
\frac{C_2}{\sqrt{n}} \biggr)\frac{t}{2\sqrt{n}\mu} \biggr]^n\nonumber\\
&= \biggl(1-\frac{\alpha^2C_1^2+\beta^2C_2^2}{4n}-\frac{(\alpha
^2+\beta^2)t^2}{2n\mu^2}-\frac{\alpha^2C_1+\beta^2C_2}{2}
\frac
{t}{\mu}\frac{1}{n} \biggr)^n \to e^{-\frac{\alpha^2{\rm Var}X+\beta
^2{\rm Var}Y}{2}}.
\end{align}
This concludes the proof of the Theorem 1.
\begin{rem}
An interesting question is
what are the implications of the assumption that $\varDelta_1^2$ and
$\varDelta_2^2$ are
negligible with respect to $\varDelta_1$ and $\varDelta_2$, Let us apply
the following formula (\cite{GR}, formula 6.616, page 710)
\begin{equation}
\int_0^{\infty}e^{-\alpha x}J_0
\bigl(\beta\sqrt{x^2+2\gamma x}\bigr)\,{\rm d}x=\frac{1}{\sqrt{\alpha^2+\beta^2}}e^{\gamma(\alpha-\sqrt
{\alpha^2+\beta^2})}.
\label{14}
\end{equation}
By using (\ref{14}), the characteristic function $\varphi(\alpha
,\beta)$ can be written as
\begin{equation}
\varphi(\alpha,\beta)\approx\frac{\mu}{\sqrt{\mu^2+\alpha
^2+\beta^2}}e^{\frac{\alpha^2\varDelta_1+\beta^2\varDelta_2}{\alpha
^2+\beta^2} (\mu-\sqrt{\mu^2+\alpha^2+\beta^2} )}.
\end{equation}
Therefore, the characteristic function of $(\hX_n,\hY_n)$, in view of
the assumptions on $\mu^{(n)}$ and $\varDelta_i^{(n)}, i=1,2$, becomes
\begin{align}
\bigl[\varphi_n(\alpha,\beta)\bigr]^n&\approx \biggl[1-\frac{\alpha
^2+\beta^2}{2(\mu^{(n)})^2}+\frac{\alpha^2\varDelta_1+\beta^2\varDelta
_2}{\alpha^2+\beta^2} \biggl(
\mu^{(n)}-\mu^{(n)}\sqrt{1+\frac
{\alpha^2+\beta^2}{(\mu^{(n)})^2}} \biggr)
\biggr]^n\nonumber\\
&\approx \biggl[1-\frac{\alpha^2+\beta^2}{2(\mu^{(n)})^2}-\frac
{\alpha^2\varDelta_1+\beta^2\varDelta_2}{\alpha^2+\beta^2}\frac{\alpha
^2+\beta^2}{2\mu^{(n)}}
\biggr]^n\nonumber\\
&\quad \to e^{-\frac{\alpha^2}{2} (\frac{t^2}{\mu^2}+\frac{C_1t}{\mu
} )-\frac{\beta^2}{2} (\frac{t^2}{\mu^2}+\frac{C_2t}{\mu
} )}. \label{15}
\end{align}
In the steps above we %considered that
made use of $\frac{\mu^{(n)}}{\sqrt{(\mu^{(n)})^2+\alpha^2+\beta
^2}}\approx1-\frac{\alpha^2+\beta^2}{2(\mu^{(n)})^2}$ as $n\to
\infty$ in view of
(\ref{20c}), and used assumption (i) afterwards.
In contrast to (\ref{5}) the terms $\frac{C_i^2}{2}$ are missing
%in the
from the limiting expression (\ref{15}).
We conclude that the approximation considered above leads to the
linearization of the limiting variances with respect to $C_1$ and $C_2$.
\end{rem}

\section{Proof of Theorem 2}\label{Theorem 2}

In the case of jump lengths with a folded Cauchy distribution (\ref{6a})
the CLT is not applicable. Again, our goal is to
%work out
calculate the limiting characteristic function %by
keeping the relevant terms in the asymptotic expansion.
We omit the lower index in $a_n, \varphi_n$ whenever it is possible.
The equalities (\ref{6}) and (\ref{7}) imply
\begin{align}
\varphi(\alpha,\beta)&=J_0 (\sqrt{\alpha^2\varDelta_1^2+\beta^2\varDelta_2^2} )\frac{2a}{\pi}\int_0^{\infty}\frac{J_0(r\sqrt{\alpha^2+\beta^2} )}{r^2+a^2}\,{\rm d}r\nonumber\\
&\quad -2J_1 (\sqrt{\alpha^2\varDelta_1^2+\beta^2\varDelta_2^2} )\frac{\alpha^2\varDelta_1+\beta^2\varDelta_2}{\sqrt{\alpha^2+\beta^2}\sqrt{\alpha^2\varDelta_1^2+\beta^2\varDelta_2^2}}\nonumber\\
&\quad \times\frac{2a}{\pi}\int_0^{\infty}\frac{J_1 (r\sqrt{\alpha^2+\beta^2})}{r^2+a^2}{\rm d}r+o(n^{-1}).
\label{20a}
\end{align}
As in Theorem 1 all the terms of the asymptotic expansion (\ref{20a})
with $k\geq2$ contain the
%factor proportional
multiplyer of
$(\alpha^2\varDelta_1^2+\beta^2\varDelta_2^2)^{k/2}$ and the remaining
sum is $o(n^{-1})$ uniformly over $(\alpha,\beta)\in U$.
Indeed, in view of (\ref{12a}) for $k\geq2$ we have
\[
\bigg|J_k \Bigl(\sqrt{\alpha^2\varDelta_1^2+
\beta^2\varDelta_2^2} \Bigr)\bigg|< C
n^{-k/2}.
\]
The module of the term $\cos(k\phi)$ in (\ref{7}) is estimated by 1.
Let $c=\sqrt{\alpha^2+\beta^2}$ and $\delta$ be the Kronecker
delta-function. Note that for any $\kappa>0$
\begin{align}
\frac{2a}{\pi}\int_0^{\infty}\frac{J_k(cr)}{r^2+a^2}\,{\rm d}r &<C\biggr(n\int_0^{n^{-2/k-\kappa}}r^{k/2}\,{\rd} r+n^{-1}\int_{n^{-2/k-\kappa}}^{\infty}\frac{J_k(cr)}{r^2}\,{\rd} r \biggl)\nonumber\\
&<C  \Bigr(n^{-\kappa-\kappa k/2}+\frac{1+\kappa}{n}\ln(n)\delta_{k=2}+n^{-1-\kappa k/2}\delta_{k>2} \Bigl),
\end{align}
and %the
a similar bound holds from below. Hence, the series is absolutely
convergent. For these reasons it remains to
%control
consider the first two terms in the expansion (\ref{20a}).

For this aim note that
\begin{equation}
\int_0^{\infty}\frac{J_1(cr)}{r^2+a^2}\,{\rd}r=
\frac{c}{a^2c^2} \Biggl[\int_0^{\infty}J_1(u)\,{
\rd}u-\int_0^{\infty}\frac
{J_1(u)u^2}{u^2+(ac)^2}\,{\rd}u \Biggr],
\label{22b}
\end{equation}
and $\int_0^{\infty}J_1(u)\,{\rd}u=-\int_0^{\infty}J'_0(u)\,{\rd
}u=1$. By differentiating
w.r.t. the parameter the formula 6.532.4 from \cite{GR} one gets
\begin{equation}
\int_0^{\infty}\frac{J_1(u)u^2}{u^2+(ac)^2}\,{\rd}u=acK_1(ac), \label{23b}
\end{equation}
where $K_1$ stands for the Macdonald function \cite{AS,W}. So, we
%obtained
obtain that
\begin{equation}
I_1=\int_0^{\infty}\frac{J_1(r\sqrt{\alpha^2+\beta
^2})}{r^2+a^2}{
\rm d}r=\frac{1}{a^2(\sqrt{\alpha^2+\beta^2)}} \bigl(1-a\sqrt{\alpha^2+
\beta^2}K_1\bigl(a\sqrt{\alpha^2+
\beta^2}\bigr) \bigr). \label{24b}
\end{equation}
If $a=a_n=\frac{\pi b}{2n}$ then for large $n$ we have $K_1(ac)=\frac
{1}{a_nc}+\frac{a_nc}{4}(2\gamma-1)+o(n^{-1})$
where $\gamma=0.57721566$ stands for the Euler--Mascheroni constant and
\[
I_1=\frac{1}{a_n^2c} \biggl[1-a_nc(
\frac{1}{a_nc}+\frac
{a_nc}{4}(2\gamma-1)+o\bigl(n^{-1}\bigr)
\biggr]=-\frac{2\gamma-1}{4}c+O\bigl(n^{-1}\bigr).
\]
As a result we obtain that the second term in (\ref{20a}) is
$O(n^{-3/2})$ and does not contribute asymptotically.

Alternatively, according to \cite{GR}, formula 6.532.1 for non-integer
$\nu$
\begin{equation}
I_{\nu}(a)=\int_0^{\infty}
\frac{J_{\nu}(x)}{x^2+a^2}\,{\rm d}x=\frac
{\pi}{a\sin(\pi\nu)} \bigl({\bar J}_{\nu}(a)-J_{\nu}(a)
\bigr), \label{25a}
\end{equation}
where ${\bar J}_{\nu}(a)$ stands for the Anger function which is a
solution of the inhomogeneous
Bessel equation $Ly=(x-\nu)\frac{\sin\pi x}{\pi}$ \cite{AS}. By
definition the Anger functions always coincide with the Bessel
functions for
the integer values of $\nu$. The following identity is well-known
\cite{DCLMT}
\begin{equation}
{\bar J}_{\nu}(x)=\frac{\sin\pi\nu}{\pi}\sum\limits
_{l=-\infty
}^{\infty}(-1)^l
\frac{J_l(x)}{\nu-l}. \label{26b}
\end{equation}
So, we use l'H\^opital's rule to evaluate the integral $I_1(a)=\lim_{\nu\to1}I_{\nu}(a)$ in (\ref{25a}) and obtain
\[
\lim_{n\to\infty}\lim_{\nu\to1}\frac{\pi}{a_n\sin(\pi\nu
)}
\bigl({\bar J}_{\nu}(a_n)-J_{\nu}(a_n)
\bigr)=0.
\]

Anyway, we need to evaluate the first term in expansion (\ref{20a}).
According to \cite{GR}, formula 6.532.6
\[
\int_0^{\infty}\frac{J_0(bx)}{x^2+a^2}\,{\rm d}x=
\frac{\pi}{2a} \bigl(I_0(ab)-L_0(ab) \bigr),
\]
where $I_0$ stands for the modified Bessel function and $L_0$ is the
modified Struve function.
Let us remind that the modified Struve function satisfies
%an
the inhomogeneous Bessel equation \cite{AS,W}
\begin{equation}
Ly=x^2\frac{d^2y}{dx^2}+x\frac{dy}{dx}-\bigl(x^2+
\nu^2\bigr)y=\frac
{4(x/2)^{\nu+1}}{\sqrt{\pi}\varGamma(\nu+1/2)}. \label{35}
\end{equation}
By using expansions of the modified Struve functions in a neighbourhood
of $0$ (see \cite{AS,W})
\begin{equation}
L_0(x)= \biggl(\frac{x}{2} \biggr)\sum
\limits_{k=0}^{\infty}
\frac
{1}{(\varGamma(3/2+k))^2} \biggl(\frac{x}{2} \biggr)^{2k},
\label{15a}
\end{equation}
we write two terms of the asymptotic expansion
\[
\frac{2a}{\pi}\int_0^{\infty}\frac{J_0 (r\sqrt{\alpha
^2+\beta^2} )}{r^2+a^2}{
\rm d}r\approx1+\frac{a^2(\alpha
^2+\beta^2)}{4}- \frac{a\sqrt{\alpha^2+\beta^2}}{2\varGamma(3/2)^2}+\cdots
\]

Finally, in view of the formula 6.565.3 of \cite{GR}
\begin{equation}
\int_0^{\infty}\frac{x^{\nu+1}}{(x^2+a^2)^{\nu+3/2}}J_{\nu
}(bx){
\rm d}x=\frac{b^{\nu}\sqrt{\pi}}{2^{\nu}a \varGamma(\nu+3/2)}e^{-ab}. \label{20}
\end{equation}
Applying (\ref{20}) for $\nu=0$ we obtain that the characteristic
function of the circular bivariate Cauchy law has the form $e^{-b \sqrt
{\alpha^2+\beta^2}}$.
So, in the limit we obtain the product of the characteristic functions
of the Gaussian law and the circular bivariate Cauchy distribution
\begin{equation}
\lim_{n\to\infty}\bigl[\varphi_n(\alpha,\beta)
\bigr]^n=\exp{ \biggl(-\frac
{C_1^2}{2}\alpha^2-
\frac{C_2^2}{2}\beta^2 \biggr)}\exp{ \bigl(-b\sqrt {
\alpha^2+\beta^2} \bigr)}, \label{21}
\end{equation}
and Theorem 2 is proved.
%% Proof %%
%%%%%%%%%%%%%%%%%%%%%%
%\begin{proof}
%\end{proof}

%% Acknowledgements %%
%%%%%%%%%%%%%%%%%%%%%%
\section*{Acknowledgments}
The article was prepared within the framework of the Academic Fund
Program at
the National Research University Higher School of Economics (HSE) and supported
%within
by the subsidy granted to the HSE
by the Government of the Russian Federation for
the implementation of the Global Competitiveness Programme.

%% Appendices %%
%%%%%%%%%%%%%%%%%%%%%%
%\appendix
\section{Appendix section}
It is interesting that for odd $2n+1$ the integral
$\int_0^{\infty}\frac{J_{2n+1}(cr)}{r^2+a^2}\,{\rm d}r$ can be
expressed in terms of the Macdonald function $K_{2n+1}(ac)$ and for
even $2n$ the integral $\int_0^{\infty}\frac
{J_{2n}(cr)}{r^2+a^2}\,{\rm d}r$ is a linear combination of the modified
Bessel function $I_{2n}(ac)$ and the modified Struve functions
$L_{k}(ac), k\leq2n$. Let us remind that the Macdonald function
$K_{\nu}(x)$
is a positive solution of the equation
\begin{equation}
Ly=x^2\frac{d^2y}{dx^2}+x\frac{dy}{dx}-\bigl(x^2+
\nu^2\bigr)y=0
\end{equation}
vanishing when $x\to\infty$ and the modified Struve function is
defined in (\ref{35}).

For odd order, Bessel functions $\int_0^{\infty}\frac
{J_{2n+1}(cr)}{r^2+a^2}\,{\rm d}r, a, c>0$ take the form
\[
2n+1=1: \frac{1}{a^2c}-\frac{K_1(ac)}{a},
\]
\[
2n+1=3: \frac{K_3(ac)}{a}+\frac{1}{a^2c}-\frac{8}{a^4c^3},
\]
\[
2n+1=5: -\frac{K_5(ac)}{a}+\frac{1}{a^2c}-\frac{24}{a^4c^3}+
\frac
{384}{a^6c^5},
\]
\[
2n+1=7: \frac{K_7(ac)}{a}+\frac{1}{a^2c}-\frac{48}{a^4c^3}+
\frac
{1920}{a^6c^5}-\frac{46080}{a^8c^7},
\]
\[
2n+1=9: -\frac{K_9(ac)}{a}+\frac{1}{a^2c}-\frac{80}{a^4c^3}+
\frac
{5760}{a^6c^5}-\frac{322560}{a^8c^7}+\frac{10321920}{a^{10}c^9},
\]
\[
2n+1=11: \frac{K_{11}(ac)}{a}+\frac{1}{a^2c}-\frac
{120}{a^4c^3}+
\frac{13440}{a^6c^5}-\frac{1290240}{a^8c^7}+\frac
{92897280}{a^{10}c^9}-\frac{3715891200}{a^{12}c^{11}}.
\]
For even order, Bessel functions $\int_0^{\infty}\frac
{J_{2n}(cr)}{r^2+a^2}\,{\rm d}r, a,c >0$ take the form
\[
2n=0: \frac{\pi I_0(ac)}{2a}-\frac{\pi L_0(ac)}{2a},
\]
\[
2n=2: \frac{c}{3}-\frac{\pi I_2(ac)}{2a}+\frac{\pi L_2(ac)}{2a},
\]
\[
2n=4: \frac{c}{15}+\frac{\pi I_4(ac)}{2a}-\frac{\pi L_2(ac)}{2a}+\frac{3\pi L_3(ac)}{a^2c},
\]
\[
2n=6: \frac{a^2c^3}{105}+\frac{c}{35}-\frac{\pi I_6(ac)}{2a}+\frac{\pi L_4(ac)}{2a}-\frac{5\pi L_3(ac)}{a^2c}+\frac{40\pi L_4(ac)}{a^3c^2},
\]
\[
2n\,{=}\,8\,{:}\, \frac{a^2c^3}{63}+\frac{c}{63}+\frac{\pi I_8(ac)}{2a}-\frac{\pi L_4(ac)}{2a}+\frac{12\pi L_5(ac)}{a^2c}-\frac{84\pi L_4(ac)}{a^3c^2}+\frac{840\pi L_5(ac)}{a^4c^3}.
\]

%% Bibliography %%
%
% structpyb loaded by rokas.maliukevicius, 2017-04-04 10:21:06

%
%%%%%%%%%%%%%%%%%%%%%%
%\bibliography{bib/paper}

\begin{thebibliography}{99}
% pybtex-1.81. Style name=spr-mp, version=1.17, label_style=nolabel, sorting_style=none, cfg=spr-mp-makeid, language=english.


%b1 ###
\bibitem{AS}
\begin{bbook}
\bauthor{\bsnm{Abramowitz}, \binits{M.}},
\bauthor{\bsnm{Stegun}, \binits{I.}}:
\bbtitle{Handbook of Mathematical Functions}.
\bpublisher{Dover}
%\blocation{???}
(\byear{1972}).
\bid{mr={0208797}}
\end{bbook}
%
\OrigBibText
Abramowitz, M., Stegun, I., Handbook of Mathematical
Functions, Dover, 1972
\endOrigBibText
\bptok{structpyb}
\endbibitem

%b2 ###
\bibitem{C}
\begin{barticle}
\bauthor{\bsnm{Chandrasekhar}, \binits{S.}}:
\batitle{Stochastic problems in physics and astronomy}.
\bjtitle{Rev. Mod. Phys.}
\bvolume{15},
\bfpage{1}--\blpage{89}
(\byear{1943}).
\bid{doi={10.1103/RevModPhys.15.1}, mr={0008130}}
\end{barticle}
%
\OrigBibText
Chandrasekhar, S., Stochastic problems in physics and astronomy, {\it
Reviews of Modern Physics}, V. 15, 1943, 1-89
\endOrigBibText
\bptok{structpyb}
\endbibitem

%b3 ###
\bibitem{DK}
\begin{bchapter}
\bauthor{\bsnm{Davydov}, \binits{Y.}},
\bauthor{\bsnm{Konakov}, \binits{V.}}:
\bctitle{Random walks in non-homogeneous Poisson environment}.
In: \beditor{\bsnm{Panov}, \binits{V.}} (ed.)
\bbtitle{Modern Problems of Stochastic Analysis and Statistics -- Festschrift in honor of Valentin Konakov}.
\bpublisher{Springer}
%\blocation{???}
(\byear{2017}).
\bcomment{In press}
\end{bchapter}
%
\OrigBibText
Davydov, Y., Konakov, V., Random walks in non-homogeneous Poisson
environment. In book: Modern Problems of Stochastic Analysis and
Statistics -
Festschrift in honor of Valentin Konakov (editor: V. Panov). Springer,
2017, In press.
\endOrigBibText
\bptok{structpyb}
\endbibitem

%b4 ###
\bibitem{GR}
\begin{bbook}
\bauthor{\bsnm{Gradshtein}, \binits{I.}},
\bauthor{\bsnm{Ryzik}, \binits{I.}}:
\bbtitle{Tables of Integrals, Sums, Series and Products}.
\bpublisher{Nauka},
\blocation{Moscow}
(\byear{1971}).
\bid{mr={0052590}}
\end{bbook}
%
\OrigBibText
Gradshtein, I., Ryzik, I., Tables of Integrals, Sums, Series and
Products, Moscow: Nauka, 1971
\endOrigBibText
\bptok{structpyb}
\endbibitem

%b5 ###
\bibitem{DGO}
\begin{bchapter}
\bauthor{\bparticle{De} \bsnm{Gregorio}, \binits{A.}},
\bauthor{\bsnm{Orsingher}, \binits{E.}}:
\bctitle{Flying randomly in $\mathbf{R}^{\mathbf{d}}$ with Dirichlet displacements}.
In: \bbtitle{Stochastic Processes Appl.},
vol.~\bseriesno{122},
pp.~\bfpage{676}--\blpage{713}
(\byear{2012}).
\bid{doi={10.1016/\\j.spa.2011.10.009}, mr={2868936}}
\end{bchapter}
%
\OrigBibText
De Gregorio, A., Orsingher, E., Flying randomly in $\bR^d$ with
Dirichlet displacements, {\it Stochastic Processes
Appl.}, V. 122, 2012, 676-713.
\endOrigBibText
\bptok{structpyb}
\endbibitem

%b6 ###
\bibitem{LC}
\begin{barticle}
\bauthor{\bparticle{Le} \bsnm{Ca\"{e}r}, \binits{G.}}:
\batitle{A Pearson random walk with steps of uniform orientation and Dirichlet distributed lengths}.
\bjtitle{J. Stat. Phys.}
\bvolume{140},
\bfpage{728}--\blpage{751}
(\byear{2010}).
\bid{doi={10.1007/s10955-010-\\0015-8}}
\end{barticle}
%
\OrigBibText
Le Ca\"er, G., A Pearson random walk with steps of uniform orientation
and Dirichlet distributed
lengths, {\it Journal of Statistical Physics}, V. 140, 2010, 728-751.
\endOrigBibText
\bptok{structpyb}
\endbibitem

%b7 ###
\bibitem{DCLMT}
\begin{barticle}
\bauthor{\bsnm{Dattoli}, \binits{G.}},
\bauthor{\bsnm{Chiccoli}, \binits{C.}},
\bauthor{\bsnm{Lorenzutto}, \binits{S.}},
\bauthor{\bsnm{Maino}, \binits{G.}},
\bauthor{\bsnm{Torre}, \binits{A.}}:
\batitle{Generalized Bessel function of the Anger type and applications to physical problems}.
\bjtitle{J. Math. Anal. Appl.}
\bvolume{184},
\bfpage{201}--\blpage{221}
(\byear{1994}).
\bcomment{N. 2}.
\bid{doi={10.1006/jmaa.1994.1194}}
\end{barticle}
%
\OrigBibText
Dattoli G., Chiccoli C., Lorenzutto S., Maino G., Torre A.,
Generalized Bessel function of the Anger type and applications to
physical problems., {\it Journ. Math. Anal. Appl.}, V.184,
1994, N. 2, 201-221
\endOrigBibText
\bptok{structpyb}
\endbibitem

%b8 ###
\bibitem{F}
\begin{bbook}
\bauthor{\bsnm{Feller}, \binits{W.}}:
\bbtitle{An Introduction to Probability Theory and its Applications, v. II}.
\bpublisher{J. Wiley \& Sons}
%\blocation{???}
(\byear{1971}).
\bid{mr={0270403}}
\end{bbook}
%
\OrigBibText
Feller W., An Introduction to Probability Theory and its Applications,
v. II, J.Wiley \& Sons, 1971
\endOrigBibText
\bptok{structpyb}
\endbibitem

%b9 ###
\bibitem{Fe}
\begin{barticle}
\bauthor{\bsnm{Ferguson}, \binits{T.}}:
\batitle{A representation of the symmetric bivariate Cauchy distribution}.
\bjtitle{Ann. Math. Stat.}
\bvolume{33}(\bissue{4}),
\bfpage{1256}--\blpage{1266}
(\byear{1962}).
\bid{doi={10.1214/aoms/1177704357}}
\end{barticle}
%
\OrigBibText
Ferguson, T., A representation of the symmetric bivariate Cauchy
distribution, {\it Ann. Math. Statist.}, V. 33, 1962,
N. 4, 1256-1266
\endOrigBibText
\bptok{structpyb}
\endbibitem

%b10 ###
\bibitem{ODG}
\begin{barticle}
\bauthor{\bsnm{Orsingher}, \binits{E.}},
\bauthor{\bparticle{De} \bsnm{Gregorio}, \binits{A.}}:
\batitle{Random flights in higher spaces}.
\bjtitle{J. Theor. Probab.}
\bvolume{20}(\bissue{4}),
\bfpage{769}--\blpage{806}
(\byear{2007}).
\bid{doi={10.1007/s10959-007-0093-y}}
\end{barticle}
%
\OrigBibText
Orsingher, E., De Gregorio, A., Random flights in higher spaces, {\it
J.Theoretical Probability}, V. 20, 2007, N.4, 769-806
\endOrigBibText
\bptok{structpyb}
\endbibitem

%b11 ###
\bibitem{W}
\begin{bbook}
\bauthor{\bsnm{Watson}, \binits{G.}}:
\bbtitle{A Treatise on the Theory of Bessel Functions}.
\bpublisher{Cambridge Univ. Press}
%\blocation{???}
(\byear{1952})
\end{bbook}
%
\OrigBibText
Watson, G., A Treatise on the Theory of Bessel Functions, Cambridge
Univ. Press, 1952
\endOrigBibText
\bptok{structpyb}
\endbibitem

\end{thebibliography}
\end{document}